\documentclass[a5paper]{book}
\usepackage[english,russian]{babel}
\usepackage{hyperref,amsthm,amsmath,amssymb,color,multirow,graphicx}
\usepackage{inputenc}
\usepackage{wrapfig, float,array}

\newtheorem{thm}{Theorem}[section]

\newtheorem{lem}[thm]{Lemma}

\theoremstyle{definition}

\numberwithin{equation}{section}

\newtheorem*{thm*}{Òåîðåìà}
\newtheorem*{lem*}{Ëåììà}

\newcommand{\T}{\mathbb T}

\newcommand{\be}{\begin{equation}}

\newcommand{\ee}{\end{equation}}

\begin{document}

\sloppy

\begin{flushright}
\begin{tabular}{l}
{\sf Uzbek  Mathematical}\\
{\sf Journal, 2021,\ No 1, pp.\pageref{12bfirstpage}-\pageref{12blastpage}}\\
\end{tabular}
\end{flushright}

\sloppy

\begin{center}
\textbf{ \large Inverse problem for a subdiffusion equation with the Caputo derivative }\\

\textbf{O.T. Muhiddinova}
\end{center}

{\small \textbf{Abstract.}  The article investigates an inverse problem of determining the right-hand side of a subdiffusion
equation with Caputo fractional derivative whose elliptic part has the most general
form and is defined on an $N$-dimensional torus $\mathbb{T}^N$. The Fourier method is used to prove theorems on the existence and uniqueness of the classical solution of the initial-boundary value problem and on the unique reconstruction of the
unknown right-hand side of the equation. Requirements for the initial function and for the additional condition are established under which the classical Fourier method can be applied to the inverse problem under
consideration.

\textbf{Keywords:}  Subdiffusion equation, the Caputo derivatives,
 Fourier method,   inverse and initial-boundary value problem.

\textbf{MSC (2010):} 35K35; 35C10. \\

\makeatletter
\renewcommand{\@evenhead}{\vbox{\thepage \hfil {\it  Muhiddinova O.T. }   \hrule }}
\renewcommand{\@oddhead}{\vbox{\hfill
{\it Inverse problem for a subdiffusion  }\hfill
\thepage \hrule}} \makeatother

\label{12bfirstpage}

\section{Introduction}

It is well known  that the interest in the
study of inverse problems for equations of mathematical physics is due to the importance of their
applications in various branches of mechanics, seismology, medical tomography, and geophysics. It
should be noted that inverse problems are well studied for differential equations of integer order.
A lot of research is also devoted to inverse problems of determining the right-hand side of subdiffusion equations (see, e.g., \cite{Zh2} - \cite{Ruz} and the bibliography therein). The paper\cite{Ruz} studied the inverse
problem of determining the right-hand side of the subdiffusion equation with Caputo fractional
derivatives whose elliptic part is an arbitrary elliptic operator with discrete spectrum and proved
the existence and uniqueness of a generalized solution by the Fourier method.  In the paper\cite{Ruz}  conditions  for  $\varphi(x)$ and $\Psi(x)$ are assumed to be defined in the operator's domain. The other papers, including \cite{Li} - \cite{Sun}, take the Caputo derivative as the fractional derivative.
  In the paper  \cite{AshurovMukhiddinova}, authors studied the inverse problem of determining the right-hand side of the subdiffusion equation with Riemann–Liouville fractional derivatives whose elliptic part is an elliptic operator of an arbitrary order defined in  bounded domain with sufficiently smooth boundary.
  
  In this article we investigate an inverse problem of determining the right-hand side of a subdiffusion
equation with Caputo fractional derivative whose elliptic part has the most general
form and is defined on an $N$-dimensional torus $\mathbb{T}^N$. The Fourier method is used to prove theorems on the existence and uniqueness of the classical solution of the initial-boundary value problem and on the unique reconstruction of the
unknown right-hand side of the equation. Requirements for the initial function and for the additional condition are established under which the classical Fourier method can be applied to the inverse problem under
consideration.

Let is move on to precise descriptions of the objects of study.

The fractional derivative in the sense of Caputo of order $0<\rho
<1$ of the function $f (t)$ defined on $[0, \infty)$ has the
form (see, for example, \cite{Pskhu}, p. 14)
$$
D_t^\rho f(t)= \frac{1}{\Gamma(1-\rho)} \int_0^t \frac{\frac{d
}{d\tau}f(\tau)d\tau}{(t-\tau)^{\rho}}, \quad t>0,
$$
provided the right-hand side exists. Here $\Gamma(\sigma)$ is
Euler's gamma function. If in this definition we interchange
differentiation and fractional integration, then we get the
definition of the Riemann-Liouville derivative:
$$
\partial_t^\rho f(t)= \frac{1}{\Gamma(1-\rho)}\frac{d}{dt}
\int_0^t \frac{f(\tau)d\tau}{(t-\tau)^{\rho}}, \quad t>0.
$$

Note that if $\rho=1$, then fractional derivatives coincides with
the ordinary classical derivative of the first order:
$\partial_tf(t) = D_t f(t)= \frac{d}{dt} f(t)$.

Let $\mathbb{T}^N$ be $N$ -dimensional torus: $\mathbb{T}^N=(-\pi,
\pi]^N$, $N\geq 1$ and
$A(D)=\sum\limits_{|\alpha|=m} a_\alpha D^\alpha$ be a homogeneous elliptic symmetric positive
differential operator with constant coefficients.

\textbf{Problem: }Find the functions $\{u(x,t), f(x)\}$ from the class: A function $f(x)\in C({T}^N)$  and a function $u(x,t)\in C({T}^N\times
[0,T])$ with the properties $D_t^\rho u(x,t), A(x,D)u(x,t)\in
C({T}^N\times (0,T])$ that satisfies the following differential equation  ($0< \rho\leq 1$)
\begin{equation}\label{eqt}
D_t^\rho u(x,t) + A(D)u(x,t) = f(x), \quad x\in \mathbb{T}^N, \quad t>0;
\end{equation}
and the initial
\begin{equation}\label{int}
 u(x,0) =\varphi(x), \quad x\in \mathbb{T}^N;
\end{equation}
 and additional conditions
\begin{equation}\label{ext}
u(x, T) = \Psi(x), \quad x\in \mathbb{T}^N.
\end{equation}

Instead of the boundary conditions we will require that the desired functions $u(x, t)$  and $f(x)$
are  $2\pi$-periodic in each variable $x_j$. We will also assume that the given functions $\varphi(x)$ and $\Psi(x)$
are $2\pi$-periodic in each variable $x_j$ as well.

\section{Main result}

To formulate the main result of this paper we recall some
definitions and statements.
 By  $A$ we denote the operator $A(D)$ defined on $2\pi$-periodic functions in $C^m(\mathbb{R}^N)$. The closure $\hat{A}$
of such an operator in $L_2(\mathbb{T}^N)$ is self-adjoint. The operator $\hat{A}$
 has a complete in $L_2(\mathbb{T}^N)$-orthonormal
system of eigenfunctions $\{(2\pi)^{-N/2}e^{inx}\}$ corresponding to the eigenvalues $A(n)$, $n\in (\mathbb{Z}^N)$.   Therefore, by virtue of the  von
Neumann theorem, for each $\tau\geq 0$ the operator
 $\hat{A}^\tau$ acts by the rule $\hat{A}^\tau
g(x)=\sum\limits_{n\in
\mathbb{Z}^N} A^\tau(n) g_n e^{inx}$,

where the $g_n$ are the Fourier coefficients of the function $g\in L_2(\mathbb{T}^N)$  in the trigonometric system
$\{(2\pi)^{-N/2}e^{inx}\}$. The domain of this operator is determined from the condition
 $\hat{A}^\tau
g(x)\in L_2(\mathbb{T}^N)$ and has the form
\begin{equation}\label{N}
D(\hat{A}^\tau)=\{g\in L_2(\mathbb{T}^N): \sum\limits_{n\in
\mathbb{Z}^N} A^{2\tau}(n) |g_n|^2 < \infty\}.
\end{equation}

To define the domain of the operator
$\hat{A}^\tau$ in terms of Sobolev spaces, recall the definition of these
spaces (see, e.g., \cite{AAP}): a function $g\in L_2(\mathbb{T}^N)$  is said to belong to the Sobolev space $L_2^a(\mathbb{T}^N)$  with
a real number $a> 0$,
 if the norm
\begin{equation}\label{T}
||g||^2_{L_2^a(\mathbb{T}^N)}=\sum\limits_{n\in
Z^N}(1+|n|^2)^a|g_n|^2.
\end{equation}
is finite. In the case of noninteger $a$ this space is also referred to as the Liouville space.

It can readily be verified that there exist constants $c_1$ and $c_2$  such that one has the estimates
$$
c_1(1+|n|^2)^{\tau m} \leq 1+A^{2\tau}(n)\leq c_2 (1+|n|^2)^{\tau
m}.
$$
Consequently, comparing the expressions (\ref{N}) and (\ref{T}),
we see that $ D(\hat{A}^\tau)= L_2^{\tau m}(\mathbb{T}^N)$.

First, let us construct a formal solution of the inverse problem (\ref{eqt}) - (\ref{ext}).
\ \ \ According to the Fourier method, we look for a solution of Problem (\ref{eqt}) - (\ref{ext}) in the form of the following series:
\begin{equation}\label{se1}
u(x,t) = \sum\limits_{n\in \mathbb{Z}^N} T_n(t) e^{inx},
\end{equation}
here the functions $T_n(t)$ are solutions of the Cauchy type problem
\begin{equation}\label{ca}
D_t^\rho T_n + A(n) T_n =f_n(t),
T_n(0)=\varphi_n.
\end{equation}
Here and everywhere below we define by $g_n$  the Fourier coefficients of the function  $g(x)$.
Note that the unknowns in problem (\ref{ca}) are also the coefficients $f_n$. To determine these numbers,
from  (\ref{ext}) we obtain the additional condition
\begin{equation}\label{ex1}
T_n(T)=\Psi_n.
\end{equation}
The unique solution of problem (\ref{ca}) has the form (see, for
example, \cite{Pskhu}, p. 16 or \cite{Gor}, p. 173)
\begin{equation}\label{Q}
T_n(t)=\varphi_n  E_{\rho, \ 1} (-A(n)
t^\rho)+ f_n\int\limits_0^t\xi^{\rho-1} E_{\rho,
\rho}(-A(n)\xi^\rho) d\xi,
\end{equation}

where $E_{\rho, \mu}$ - is the Mittag-Leffler function
$$
E_{\rho, \mu}(t)= \sum\limits_{k=0}^\infty \frac{t^k}{\Gamma(\rho
k+\mu)}
$$
and  $\Gamma (k)$ - is the Euler's gamma function.

For the integral one has the formula (see, for example. \cite{Gor}, formula (4.4.4))
\begin{equation}\label{ML0}
\int\limits_0^t \xi^{\rho-1} E_{\rho, \rho}(-\ A(n)\xi^\rho)
d\xi=t^\rho\, E_{\rho, \rho+1}(-\ A(n) t^\rho).
\end{equation}
and we have
\begin{equation}\label{Tk}
T_n(t)=\varphi_n  E_{\rho, \ 1} (-A(n)
t^\rho)+ f_nt^\rho\, E_{\rho, \rho+1}(-\ A(n) t^\rho).
\end{equation}
To determine the coefficient $f_n$, by virtue of the additional condition
(\ref{ex1}) we obtain the equation
$$
T_n(T)=\varphi_n  E_{\rho, \ 1} (-A(n)
T^\rho)+ f_nT^\rho\, E_{\rho, \rho+1}(-\ A(n) T^\rho)=\Psi_n.
$$
 Consequently,
\begin{equation}\label{fk}
f_n=\frac{\Psi_n}{T^\rho E_{\rho,
\rho+1}(-A(n)T^\rho)}-\frac{\varphi_n E_{\rho, \ 1}(-A(n)
T^\rho)}{T^\rho E_{\rho, \rho+1}(-A(n) T^\rho)}.
\end{equation}
Thus, we have constructed the formal series (see (\ref{se1}) and (\ref{Tk}))
\begin{equation}\label{u}
u(x,t)=\sum\limits_{n\in\mathbb{Z}^N}\varphi_n
E_{\rho, \ 1} (-A(n) t^\rho)e^{inx}+\sum\limits_{n\in\mathbb{Z}^N} f_nt^{\rho} E_{\rho,
\rho+1}(-A(n)t^\rho) e^{inx},
\end{equation}

for the solution of the inverse problem, where the numbers $f_n$ are determined using formula (\ref{fk}) and
\begin{equation}\label{f}
f(x)=\sum\limits_{n\in\mathbb{Z}^N}\frac{\Psi_n}{T^\rho E_{\rho,
\rho+1}(-A(n)T^\rho)}e^{inx}-\sum\limits_{n\in\mathbb{Z}^N} \frac{\varphi_n E_{\rho, \ 1}(-A(n)
T^\rho)}{T^\rho E_{\rho, \rho+1}(-A(n) T^\rho)}e^{inx}.
\end{equation}
Obviously, the function (\ref{u}) and (\ref{f}) are formal solutions

Now let us state the following theorem on the existence of a solution of the inverse problem on
the torus.

\begin{thm}\label{TN} Let $\varphi\in L^{\tau}_2(\mathbb{T}^N)$ and
$\Psi\in L^{\tau+1}_2(\mathbb{T}^N)$, where $\tau > \frac{N}{2}$.
 Then there exists a unique
periodic solution  $\{u(x,t),
f(x)\}$  of the inverse problem (\ref{eqt}) -  (\ref{ext})  which can be represented in the
form of the series (\ref{u}) and (\ref{f}).

\end{thm}

\section{Uniqueness}

\begin{thm}\label{Ed} (on uniqueness).  Let the function $\varphi(x)$ and $\Psi(x)$
be continuous in  ${\mathbb{T}^N}$. Then there can exist only one classical solution $\{u(x,t), f(x)\}$
of the inverse problem (\ref{eqt})- (\ref{ext}).
\end{thm}

\noindent \textbf{Proof. }

Assume that under the assumptions in the theorem there are two pair of solutions, $\{u_1(x,t), f_1(x)\}$ and $\{u_2(x,t), f_2(x)\}$. Let us prove that $u(x,t)=u_1(x,t)-u_2(x,t)\equiv 0$ and $f(x)=f_1(x)-f_2(x)\equiv
0$.  Since the problem under consideration is linear, we have the following problem for determining the functions
$u(x,t)$ and $f(x)$ :
\begin{equation}\label{eq1}
D_t^\rho u(x,t) + A(D)u(x,t) = f(x), \quad x\in T^{N}, \quad t>0;
\end{equation}
\begin{equation}\label{in2}
u(x,0)=0, \quad
x\in T^{N},
\end{equation}
\begin{equation}\label{ex0}
u(x, T)=0, \quad x\in T^{N}.
\end{equation}

Here, instead of boundary conditions, we require the functions we
are looking for to be periodic with respect to each variable.

Let $u\left(x,t\right)$ be a solution of problem
\eqref{eq1} - \eqref{ex0}.
Consider the following function
\begin{equation} \label{GrindEQ__3_3_}
\ w _{n} \left(t\right)=(2\pi) ^{-{\raise0.7ex\hbox{$ N
$}\!\mathord{\left/ {\vphantom {N 2}} \right.
\kern-\nulldelimiterspace}\!\lower0.7ex\hbox{$ 2 $}} } \int
_{T^{N} }u\left(x,t\right) \cdot e^{inx} dx.
\end{equation}

If we differentiate this equality under the integral and take into account equation \eqref{eq1}, we obtain the following equality
\[
D_t^\rho w_n(t)=(2\pi) ^{-{\raise0.7ex\hbox{$ N $}\!\mathord{\left/
{\vphantom {N 2}} \right.
\kern-\nulldelimiterspace}\!\lower0.7ex\hbox{$ 2 $}} } \int
_{T^{N} }D_t^\rho u \left(x,t\right) \cdot e^{inx} dx=\]
\[=-(2\pi)
^{-{\raise0.7ex\hbox{$ N $}\!\mathord{\left/ {\vphantom {N 2}}
\right. \kern-\nulldelimiterspace}\!\lower0.7ex\hbox{$ 2 $}} }
\int _{T^{N} }A\left(D\right)u\left(x,t\right)\cdot  e^{inx} dx + (2\pi) ^{-{\raise0.7ex\hbox{$ N
$}\!\mathord{\left/ {\vphantom {N 2}} \right.
\kern-\nulldelimiterspace}\!\lower0.7ex\hbox{$ 2 $}} } \int
_{T^{N} }f(x) \cdot e^{inx} dx,
t>0.
\]

Integrating here by parts and taking into account that the operator $\hat{A}$ is self-adjoint, we obtain the following equality
\[
D_t^\rho w_n(t)=-(2\pi) ^{-{\raise0.7ex\hbox{$ N $}
\!\mathord{\left/ {\vphantom {N 2}} \right.
\kern-\nulldelimiterspace}\!\lower0.7ex\hbox{$ 2 $}} } \int
_{T^{N} }u\left(x,t\right)A\left(D\right) \cdot e^{inx}
dx+f_n\]
\[=-A(n)\cdot (2\pi) ^{-{\raise0.7ex\hbox{$ N $} \!\mathord{\left/
{\vphantom {N 2}} \right.\kern-\nulldelimiterspace}\!\lower0.7ex\hbox{$ 2 $}} } \int
_{T^{N} }u\left(x,t\right) \cdot e^{inx} dx+f_n;\]
\[D_t^\rho w_n(t)=-A(n)\ w _{n} \left(t\right)+f_n,t>0.
\]
Now, taking into account initial condition \eqref{in2}, we conclude that the function $\ w_{n}\left(t\right)$ is the solution of the following Cauchy problem:
\[D_t^\rho w_n(t)+A(n)\ w _{n} \left(t\right)-f_n=0,     t>0;    \ w_{n} \left(0\right)=0.\]

The solution of this problem has the form (see, for
example,  \cite{Pskhu} p. 16 or \cite{Gor}, p. 173 )
$$
w_n(t)=f_n\int\limits_0^t \xi^{\rho-1} E_{\rho,
\rho}(-\ A(n)\xi^\rho) d\xi,
$$
 For the integral one has the formula (see, for example. \cite{Gor}, formula (4.4.4))
\begin{equation}\label{ML}
\int\limits_0^t \xi^{\rho-1} E_{\rho, \rho}(-\ A(n)\xi^\rho)
d\xi=t^\rho\, E_{\rho, \rho+1}(-\ A(n) t^\rho).
\end{equation}
Therefore, using the additional condition (\ref{ex0}), we obtain $w_n(T)=f_n\,T^\rho\, E_{\rho, \rho+1}(-\ A(n) T^\rho)=0$
and since $0<E_{\rho, \rho+1}(-\ A(n) T^\rho)\leq 1$, we conclude that
$f_n=0$ for all $n$. In turn, this implies that $w_n(t)\equiv
0$ for all $n$. Consequently, by virtue of the completeness of the system of eigenfunctions $\{(2\pi)^{-N/2}e^{inx}\}$, we have  $f(x)=0$ and $u(x,t) = 0$ for all $x\in
\mathbb{T}^N$ and $t>0$.  Thus, we have proved the uniqueness of
solution of the inverse problem (\ref{eqt})- (\ref{ext}) .

\section{Existence}

\noindent The proof of existence is based on the following lemma, which is, in fact, corollary of the classical Sobolev embedding theorem (see M.A. Krasnosel'skii et.al., \cite{st3}, p. 453).

\begin{lem}\ Let $\tau >\frac{|\alpha|}{m}+\frac{N}{2m}$
\textit{. Then for any $\left|\alpha \right|\le m,$ the operator
}$D^{\alpha } \left(\hat{A}+{\rm I} \right)^{-\tau } $\textit{
maps (completely) continuously from} \textit{$L_{2} \left(T^{N}
\right)$ into $C\left(T^{N} \right)$, moreover the following
estimate holds true:}
\begin{equation} \label{GrindEQ__4_1_}
\left\| D^{\alpha } \left(\hat{A}+{\rm I} \right)^{-\tau }
g\right\| _{C\left(T^{N} \right)} \le C\left\| g\right\| _{L_{2}
\left(T^{N} \right)}
\end{equation}
\end{lem}
\textit{where we denote by ${\rm I} $ the unit operator in} $L_{2}
\left(T^{N} \right)$.

\noindent \textbf{Proof. }Since the embedding theorem one has

$$\left\| D^{\alpha } \left(\hat{A}+{\rm I}
\right)^{-\tau } g\right\| _{C\left(T^{N} \right)} \le C\left\|
D^{\alpha } \left(\hat{A}+{\rm I} \right)^{-\tau } g\right\|
_{L_{2}^{a} \left(T^{N} \right)}, $$ 
for any \textit{$a>\frac{N}{2}
$}. Therefore it is sufficient to prove the inequality 
$$\left\|
D^{\alpha } \left(\hat{A}+{\rm I} \right)^{-\tau } g\right\|
_{L_{2}^{a} \left(T^{N} \right)} \le C\left\| g\right\| _{L_{2}
\left(T^{N} \right)}. $$

\noindent But this is a consequence of the estimate
\[\sum _{n\in Z^{n} }\left|g_{n} \right|^{2} \left|n\right| ^{2\left|\alpha
\right|} (1+A\left(n\right))^{-2\tau } \left(1+\left|n\right|^{2} \right)^{a} \le C\sum _{n\in Z^{n} }\left|g_{n} \right|^{2}  ,\]
that is valid for \textit{$\frac{N}{2} <a\le \tau m-\left|\alpha
\right|$}.

\noindent Now we show how the existence of the solution of problem (\ref{eqt})- (\ref{ext}) can be derived from this lemma.

Let the assumptions of  Theorem 2.1 be satisfied. First, let us prove the absolute and uniform
convergence of the series  (\ref{f}) for the function $f$. Set
\begin{equation}\label{f1}
f^1_j(x)=\sum _{n\in\mathbb{Z}^N}\frac{\Psi_n}{T^\rho E_{\rho,
\rho+1}(-A(n)T^\rho)}e^{inx},
\end{equation}
\begin{equation}\label{f2}
f^2_j(x)=\sum _{n\in\mathbb{Z}^N} \frac{\varphi_n E_{\rho, \ 1}(-A(n)
T^\rho)}{T^\rho E_{\rho, \rho+1}(-A(n) T^\rho)}e^{inx}.
\end{equation}

The condition of the theorem implies the convergence of the series
$$
\sum\limits_{n\in\mathbb{Z}^N}\ (A(n)+1)^{2(\tau+1)} |\Psi_n|^2 \leq
C_\Psi<\infty,
$$

 for some $\tau> \frac{N}{2m}$and the convergence of the series
$$
\sum\limits_{n\in\mathbb{Z}^N}\ (A(n)+1)^{2\tau} |\varphi_n|^2 \leq
C_\varphi<\infty.
$$
for the function $\varphi$.

Since $(\hat{A}+I )^{-\tau } e^{inx} =(A(n)+1)^{-\tau }e^{inx}$, for $f^1_j$ we have
$$
f^1_j(x)=(\hat{A}+I )^{-\tau }\sum _{n\in\mathbb{Z}^N}\frac{\Psi_n\, \ (A(n)+1)^\tau}{T^\rho E_{\rho,
\rho+1}(-A(n)T^\rho)}e^{inx}.
$$
Then, applying the lemma  we obtain
$$
|| f^1_j||_{C\left(T^{N} \right)} = ||(\hat{A}+I )^{-\tau }\sum _{n\in\mathbb{Z}^N}\frac{\Psi_n\, \ (A(n)+1)^\tau}{T^\rho E_{\rho,
\rho+1}(-A(n)T^\rho)}e^{inx}||_{C\left(T^{N} \right)}\leq
$$
\begin{equation}\label{f1L}
\leq C ||\sum _{n\in\mathbb{Z}^N}\frac{\Psi_n\, \ (A(n)+1)^\tau}{T^\rho E_{\rho,
\rho+1}(-A(n)T^\rho)}e^{inx}||_{L_2\left(T^{N} \right)}.
\end{equation}
Based on this, by virtue of the orthonormality of the system $\{(2\pi)^{-N/2}e^{inx}\}$,  we have
$$
|| f^1_j||^2_{C\left(T^{N} \right)} \leq C \sum _{n\in\mathbb{Z}^N}
\bigg|\frac{\Psi_n\, \ (A(n)+1)^\tau}{T^\rho E_{\rho,
\rho+1}(-A(n)T^\rho)}\bigg|^2.
$$

Let us give estimates for the Mittag-Leffler function required in the sequel. For large $t$  one has
the asymptotic estimate (see e.g.,
\cite{Pskhu}, p. 13, recall that $0<\rho\leq 1$)
\begin{equation}\label{ML1}
E_{\rho, \rho+1}(-t)= \frac{1}{t}\big(1+O(\frac{1}{t})\big), \quad
t>1,
\end{equation}
and for arbitrary complex numbers $\mu$ we have the estimate
\begin{equation}\label{ML2}
0< |E_{\rho, \mu}(-t)|\leq \frac{C}{1+t}, \quad t>0.
\end{equation}

Let us get back to estimating the norm of the function $f_j^1$. By virtue of the estimate
(\ref{ML1}), we have
\begin{equation}\label{f1C}
|| f^1_j||^2_{C\left(T^{N} \right)} \leq C \sum\limits_{n\in\mathbb{Z}^N}\ (A(n)+1)^{2(\tau+1)} |\Psi_n|^2 \leq
C_\Psi<\infty.
\end{equation}
Applying the same argument as above to the function$f_j^2$,  we arrive at the inequality
\begin{equation}\label{f2L}
|| f^2_j||_{C\left(T^{N} \right)} \leq C ||\sum _{n\in\mathbb{Z}^N} \frac{\varphi_n \ (A(n)+1)^\tau E_{\rho, \ 1}(-A(n)
T^\rho)}{T^\rho E_{\rho, \rho+1}(-A(n) T^\rho)}e^{inx}||_{L_2\left(T^{N} \right)}.
\end{equation}
Hence, by virtue of the orthonormality of the system $\{(2\pi)^{-N/2}e^{inx}\}$,  we have
$$
|| f^2_j||^2_{C\left(T^{N} \right)} \leq C \sum _{n\in\mathbb{Z}^N}
\bigg|\frac{\varphi_n \ (A(n)+1)^\tau E_{\rho, \ 1}(-A(n)
T^\rho)}{T^\rho E_{\rho, \rho+1}(-A(n) T^\rho)}\bigg|^2.
$$
Further, applying the estimates (\ref{ML1}) and (\ref{ML2}), we obtain
\begin{equation}\label{f2C}
|| f^2_j||^2_{C\left(T^{N} \right)} \leq \frac{C}{T^\rho}
\sum_{n\in\mathbb{Z}^N} \ (A(n)+1)^{2\tau} |\varphi_n|^2 \leq
\frac{C}{T^\rho} C_\varphi<\infty.
\end{equation}

The estimates (\ref{f1C}) and (\ref{f2C}) imply the uniform convergence of series (\ref{f}),defining the function  $f(x)$. On the other hand, the sums in (\ref{f1L}) and (\ref{f2L}) converge under any permutation of their
terms, because these terms are mutually orthogonal. This implies the absolute convergence of 
series (\ref{f}).

Let us proceed to considering the series in (\ref{u}) for the function
$u$. To this end, using the definition
of the coefficients $f_n$,  we study the following three sums:

\begin{equation}\label{u1}
u^1_j(x,t)=\sum _{n\in\mathbb{Z}^N} \varphi_n  E_{\rho, \ 1}
(-\ A(n) t^\rho)e^{inx},
\end{equation}
\begin{equation}\label{u2}
u^2_j(x,t)=\sum _{n\in\mathbb{Z}^N} \frac{\Psi_n t^{\rho} E_{\rho,
\rho+1} (-\ A(n) t^\rho)}{T^\rho E_{\rho,
\rho+1}(-\ A(n)T^\rho)} e^{inx},
\end{equation}
\begin{equation}\label{u3}
u^3_j(x,t)=\sum _{n\in\mathbb{Z}^N} \frac{\varphi_n E_{\rho,
\ 1}(-\ A(n) T^\rho)}{T^\rho E_{\rho, \rho+1}(-\ A(n)
T^\rho)}\,t^{\rho} E_{\rho, \rho+1} (-\ A(n) t^\rho) e^{inx}.
\end{equation}

Let $|\alpha|= m$ and $\delta>0$. We need to show that each of these sums converges uniformly and
absolutely in the domain $\mathbb\T^{N}\times [\delta, T]$ after the operators $D^\alpha$and $D^\rho_t$
are term-by-term applied to the sum.

We have
$$
u^1_j(x,t)=(\hat{A}+I)^{-\tau-1}\sum _{n\in\mathbb{Z}^N}
\ (A(n)+1)^{\tau+1}\varphi_n E_{\rho,\ 1}(-\ A(n)
t^\rho)e^{inx}.
$$
Therefore, by the lemma, one has the estimate
$$
||D^\alpha u^1_j||_{C\left(T^{N} \right)}=||D^\alpha
(\hat{A}+I)^{-\tau-1}\sum _{n\in\mathbb{Z}^N}
\ (A(n)+1)^{\tau+1}\varphi_n E_{\rho,\ 1}(-\ A(n)
t^\rho)e^{inx}||_{C\left(T^{N} \right)}\leq
$$
\begin{equation}\label{u1L}
\leq C ||\sum _{n\in\mathbb{Z}^N}
\ (A(n)+1)^{\tau+1}\varphi_n E_{\rho,\ 1}(-\ A(n)
t^\rho)e^{inx}||_{L_2(\left(T^{N} \right)}.
\end{equation}
Using the orthonormality of the system $\{(2\pi)^{-N/2}e^{inx}\}$,  we have
$$
||D^\alpha u^1_j||^2_{C\left(T^{N} \right)}\leq C \sum _{n\in\mathbb{Z}^N}
|\ (A(n)+1)^{\tau+1}\varphi_n E_{\rho,\ 1}(-\ A(n)
t^\rho)|^2\leq
$$
$$
\leq C t^{-2\rho}\sum _{n\in\mathbb{Z}^N}
\ (A(n)+1)^{2\tau}|\varphi_n|^2\leq C \delta^{-2\rho}C_\varphi.
$$
Here the last inequality has been obtained based on the estimate (\ref{ML2}),

This implies the uniform convergence of the sum (\ref{u1}) differentiated with respect to the variables $x_n$ in the closed cylinder $(x,
t)\in\mathbb\T^{N}\times [\delta, T]$.   On the other hand, the sum (\ref{u1L}) converges
under any permutation of its terms, because these terms are mutually orthogonal. This implies the absolute convergence of the sum (\ref{u1})  differentiated with respect to the variables $x_n$.

Repeating a similar argument, we establish that the same assertions also hold true for the
sums (\ref{u2}) and (\ref{u3}),  determining the functions
$u_2(x,t)$ and $u_3(x,t)$.

Now it remains to verify the validity of term-by-term application of the operator $D_t^\rho$ to the
series  (\ref{u}) 
defining the solution $u(x,t)$. It can readily be verified that
$$
D_t^\rho\sum\limits_{n\in \mathbb{Z}^N} T_n(t) e^{inx}=
-\sum\limits_{n\in \mathbb{Z}^N} A(n) T_n(t) e^{inx}+\sum\limits_{n\in \mathbb{Z}^N} f_n(t) e^{inx}=
$$
$$
-A(x,D)(\hat{A}+I)^{-\tau-1}\sum\limits_{n\in \mathbb{Z}^N} (\ A(n)+1)^{\tau+1}
T_n(t)e^{inx}+(\hat{A}+I)^{-\tau}\sum\limits_{n\in \mathbb{Z}^N}(\ A(n)+1)^{\tau}f_k e^{inx}.
$$
The absolute and uniform convergence of the last two sums has already been proved above.

The fact that the functions (\ref{u}) and (\ref{f}) satisfy all conditions (\ref{eqt})- (\ref{ext}) of the inverse problem
is obvious. Theorem has thus been proved completely.

\begin{center}
\textbf{Acknowledgements}
\end{center}

\noindent The author convey thanks to R. R. Ashurov for
discussions of these results.

\begin{tabular}{p{9cm}}

\\Muhiddinova O.T.,\\

\\ Institute of Mathematics named after V.I. Romanovsky, the Academy of Sciences of the Republic of Uzbekistan, University street, 4, Tashkent, 100174, Uzbekistan\\

\\ e-mail: oqila1992@mail.ru\\

\end{tabular}

\end{document}